\documentclass{article}

\usepackage{amsfonts}
\usepackage{amssymb}
\usepackage{amsbsy}
\usepackage{amsmath}
\usepackage{amsthm} 
\usepackage{url}

\newtheorem{theorem}{Theorem}
\newtheorem{conjecture}[theorem]{Conjecture}
\newtheorem{definition}[theorem]{Definition}
\newtheorem{proposition}[theorem]{Proposition}

\newcommand{\M}{\mathcal{M}}
\newcommand{\N}{\mathbf{N}}
\newcommand{\Z}{\mathbf{Z}}
\newcommand{\C}{\mathbf{C}}
\newcommand{\Q}{\mathbf{Q}}
\newcommand{\kk}{\mathbf{k}}

\newcommand{\Magma}{{\sc Magma}}
\newcommand{\Sage}{{\sc Sage}}

\newcommand{\Gouvea}{Gouv\^ea}

\DeclareMathOperator{\SL}{SL}

\newcommand{\sltz}{\SL_2(\Z)}

\title{On the~$U_p$ operator acting on $p$-adic overconvergent modular forms when~$X_0(p)$ has genus~1}

\author{L. J. P. Kilford}

\begin{document}

\maketitle

\begin{abstract}
In this note we will show how to compute $U_p$ acting on spaces of overconvergent $p$-adic modular forms when $X_0(p)$ 
has genus~1. We first give a construction of Banach bases for spaces of overconvergent $p$-adic modular forms, and then 
give an algorithm to approximate both the characteristic power series of the~$U_p$ operator and eigenvectors of finite 
slope for~$U_p$, and present some explicit examples. We will also relate this to the conjectures of Clay on the slopes of overconvergent modular forms.\end{abstract}

\section{Introduction}

Let~$p$ be a prime number. Let us first recall the definition of the \emph{slope} of a normalized modular eigenform.
\begin{definition}
Let~$f$ be a normalized cuspidal modular eigenform with Fourier expansion at~$\infty$ given by~$\sum_{n=1}^\infty a_n 
q^n$ (so~$a_1=1$). The \emph{slope} of~$f$ is defined to be the $p$-valuation of~$a_p$ viewed as an element of~$\C_p$; 
we normalize the $p$-valuation of~$p$ to be~$1$. We also say that the slope of~0 is infinite. \end{definition}

There has been substantial interest in computing the action of the~$U_p$ operator acting on overconvergent modular forms in the case where~$X_0(p)$ has genus~0; in other words, when~$p \in \{2,3,5,7,13\}$; indeed, computations such as that in~\cite{kolberg} can be seen as precursors to this. One particular approach to this was outlined in~\cite{smithline} (see~\cite{loeffler} and~\cite{smithline-published} for an accessible exposition of this), and there are also explicit results on the slopes of the~$U_p$ operator in certain circumstances; see~\cite{buzzard-calegari}, \cite{buzzard-kilford}, \cite{coleman-stevens-teitelbaum}, \cite{emerton-thesis} and~\cite{kilford-5slopes}, as well as conjectures which give formulae for the slopes; see~\cite{buzzard-questions} and~\cite{clay-thesis}, for instance. In this paper, we will give a similar construction for those cases where~$X_0(p)$ has genus~1.

There have also been a number of papers which attempt to approximate overconvergent modular forms or automorphic forms, such as~\cite{atkin-obrien}, \cite{gouvea-mazur-searching} and~\cite{loeffler}. These either prove or assume the existence of an eigenform for the~$U_p$ operator, and then use an iterative process to find an approximation to the eigenform of the desired accuracy. In this paper, we will give a construction which allows us to do this.

\section{Some arithmetic geometry}
Throughout the rest of this paper, we will assume that~$p$ is a prime number such that~$X_0(p)$ has genus~1; in other 
words, that~$p \in \{11,17,19\}$. We will also only be concerned with forms with finite slope (in other words, with nonzero $U_p$ eigenvalue); 
see~\cite{coleman-stein} for an investigation of when forms of infinite slope can be approximated by forms of finite slope.

\subsection{Background material}
We first recall the definition of overconvergent $p$-adic modular forms; we recall that these are defined to give a subspace of the $p$-adic forms on which the Hecke operator~$U_p$ is compact. This will require us to summarize some background material.

Following Katz~\cite{katz}, section~2.1, we recall that, for~$C$ an elliptic curve over an $\mathbf{F}_p$-algebra~$R$, there is a mod~$p$ modular form~$A(C)$ called the \emph{Hasse invariant}, which has $q$-expansion over~$\mathbf{F}_p$ equal to~$1$.

We consider the Eisenstein series of weight~$p-1$ and tame level~1 defined over~$\mathbf{Q}$, with Fourier expansion
\[
E_{p-1}(q):=1+\frac{2(p-1)}{B_{2(p-1)}} \sum_{n = 1}^\infty \left( \sum_{0 < d | n} d^{p-1} \right) \cdot q^n.
\]
We see that~$E_{p-1}$ is a lifting of~$A(C)$ to characteristic~0, as the reduction of~$E_{p-1}$ to characteristic~$p$ 
has the same $q$-expansion as~$A(C)$, and therefore~$E_{p-1} \mod p$ and~$A(C)$ are both modular forms of level~1 and 
weight~$p-1$ defined over~$\mathbf{F}_p$, with the same $q$-expansion. Note also that if~$C$ is an elliptic curve 
defined over~$\Z_p$ then the valuation~$v_p(E_{p-1}(C))$ can be shown to be well-defined. We also note that for the 
primes~$p$ we are considering here, the Eisenstein series does not necessarily have coefficients in~$\Z$ (this is 
because the space of modular forms for~$\sltz$ of weight~$p-1$ is not necessarily 1-dimensional; this is another new 
feature that we have to deal with in this paper).

We now let~$m$ be a positive integer, and recall that the modular curve $X_0(p^m)$ is defined to be the moduli space that parametrizes pairs~$(C,P)$, where~$C$ is an elliptic curve and~$P$ is a subgroup of~$C$ of order~$p^m$. 

Using arguments exactly similar to those in~\cite{kilford-2slopes}, we define the affinoid subdomain~$Z_0(p^m)$ of~$X_0(p^m)$ to be the connected component containing the cusp~$\infty$ of the set of points~$t=(C,P)$ in~$X_0(p^m)$ which have~$v_p(E_{p-1}(t))=0$; we note that~$v_p(E_{p-1}(t))=0$ means either that the point~$t$ corresponds to an ordinary elliptic curve or that~$t$ is a cusp. 

We now define strict affinoid neighbourhoods of~$Z_0(p^m)$.

\begin{definition}[Coleman~\cite{coleman}, Section~B2]
\label{connected-component}
We think of~$X_0(p^m)$ as a rigid space over~$\mathbf{Q}_p$, and we let~$t \in X_0(p^m)(\overline{\mathbf{Q}}_p)$ be a point, corresponding either to an elliptic curve defined over a finite extension of~$\Q_p$, or to a cusp.
Let~$w$ be a rational number, such that~$0 < w < p^{2-m}/(p+1)$.

We define~$Z_0(p^m)(w)$ to be the connected component of the affinoid
\[
\left\{t \in X_0(p^m): \; v_p(E_{p-1}(t)) \le w\right\}
\]
which contains the cusp~$\infty$.
\end{definition}

\subsection{Overconvergent $p$-adic modular forms}

Having introduced this terminology, we can now define $p$-adic overconvergent modular forms.

\begin{definition}[Coleman,~\cite{coleman-overconvergent}, page~397]
Let~$w$ be a rational number, such that~$0 < w < p^{2-m}/(p+1)$. Let~$\mathcal{O}$ be the structure sheaf of~$Z_0(p^m)(w)$.
We call sections of~$\mathcal{O}$ on~$Z_0(p^m)(w)$ 
\emph{$w$-overconvergent $p$-adic modular forms of weight~$0$ and level~$\Gamma_0(p^m)$}.
If a section~$f$ of~$\mathcal{O}$ is a $w$-overconvergent modular form, then we say that~$f$ is an \emph{overconvergent $p$-adic modular form}.

Let~$K$ be a complete subfield of~$\mathbf{C}_p$, and define~$Z_0(p^m)(w)_{/K}$ to be the affinoid over~$K$ induced from~$Z_0(p^m)(w)$ by base change from~$\mathbf{Q}_5$. The space$$M_0(p^m,w;K):=\mathcal{O} (Z_0(p^m)(w)_{/K})$$ of $w$-overconvergent modular forms of weight~$0$ and level~$\Gamma_0(p^m)$ is a $K$-Banach space.

We now let~$\chi$ be a primitive Dirichlet character of conductor~$p^m$ and let~$k$ be an integer such that~$\chi(-1)=(-1)^k$. Let~$E^*_{k,\chi}$ be the normalized Eisenstein series of weight~$k$ and character~$\chi$ with nonzero constant term.

The space of $w$-overconvergent $p$-adic modular forms of weight~$k$ and character~$\chi$ is given by
\begin{equation}
\label{eq:weight-k-overconvergent-forms}
\M_{k,\chi}(p^m,w;K):=E^*_{k,\chi} \cdot \M_0(p^m,w;K).
\end{equation}
This is a Banach space over~$K$.
\end{definition}

There are Hecke operators~$U_p$ and~$T_l$ (where~$l\ne p$) acting on the space of modular forms~$\M_{k,\chi}(p^m,w;K)$; these are defined on the $q$-expansions of the overconvergent modular forms in exactly the same way as they are defined on the Fourier expansions of classical modular forms. One defines~$T_n$ for~$n$ a natural number in the usual way.

Using results of Coleman, we have the following theorem about the independence of the characteristic power series 
of~$U_p$ acting on~$\M_{k,\chi}(p^m,w;K)$, which means that we can prove results for a single choice of~$w$ and know 
that our result does not depend on this choice. This is the result which allows us to draw conclusions about the 
characteristic power series of~$U_p$ from our computations, which involve choosing one value of~$w$ and working with 
that value.
\begin{theorem}[Coleman~\cite{coleman}, Theorem~B3.2]
Let~$w$ be a real number such that~$0 < w < \min(p^{2-m}/p,1/(p+1))$, let~$k$ be an integer and let~$\chi$ be a character such that~$\chi(-1)=(-1)^k$.

The characteristic polynomial of~$U_p$ acting on $w$-overconvergent $p$-adic modular forms of weight~$k$ and character~$\chi$ is independent of the choice of~$w$.
\end{theorem}

Because we will be approximating the matrix of the~$U_p$ operator in the next section, we will find the following 
result of Serre useful (we note that in~\cite{serre} the older term 
``completely continuous'' is used). We note also 
that it is a well-known result that the~$U_p$ operator acting on overconvergent modular forms is compact; we recall 
that this was the motivation for introducing overconvergent modular forms as a subspace of the $p$-adic modular 
forms. 
\begin{theorem}[Serre~\cite{serre}, Proposition~7]
\label{serre-proposition}
Let~$M_\infty$ be a compact infinite matrix (that is, the matrix of a compact operator). If~$M_m$ is a series of finite matrices which tend to~$M_\infty$, then the finite characteristic power series~$\mathop{det}(1-tM_m)$ converge coefficientwise to~$\mathop{det}(1-tM_\infty)$, as~$m \rightarrow \infty$.
\end{theorem}
If the genus of~$X_0(p)$ is zero, then it can be shown that a Banach basis for the space of $p$-adic overconvergent modular forms of weight~0 is 
\[
\{z,z^2,\ldots,z^i,\ldots\},
\]
where~$z$ is a suitable modular function of weight~0 which vanishes at the 
cusp~$\infty$. The radius of overconvergence will depend on the choice 
of~$z$; see~\cite{kilford-2slopes} and~\cite{kilford-5slopes} for 
examples of this.

However, one cannot use the same methods when the genus of~$X_0(p)$ is greater than~0, because there are now 
several supersingular $j$-invariants modulo~$p$ (in general, there are approximately~$\lfloor p/12\rfloor$ modulo~$p$). In the genus~1 case, which we are dealing with here, there are \emph{two} 
supersingular $j$-invariants; these can be computed to be~0 and~1 when~$p=11$, 0 
and~8 when~$p=17$ and~7 and~18 when~$p=19$. This means that the overconvergent functions of weight~0 are functions on 
an annulus rather than on a disc (as in the genus~0 case). 

We will now state the main result of this paper. 
\begin{proposition}
Let~$p$ be a prime such that~$X_0(p)$ has genus~1 and let~$a$ and~$b$ be 
the distinct supersingular $j$-invariants modulo~$p$. Let~$z:=c \cdot (j-a)/(j-b)$ 
be a parameter on the closed unit disc, where~$c \in \C_p$ has normalized valuation strictly between~0 and~1.
 
If we let~$\mathcal{A}$ be the annulus $\{1/p\le|z|\le1\}$ in this disc, 
then an orthonormal basis for the functions on~$\mathcal{A}$ is 
given by
\begin{equation}
\label{eq:orthonormal-basis}
\{1,z,z^2,z^3,\ldots,z^i,\ldots\} \cup 
\{p/z,(p/z)^2,(p/z)^3,\ldots,(p/z)^i,\ldots\}.
\end{equation}
These are \emph{$p$-adic overconvergent modular forms of weight~0}; we 
call the 
space spanned by nonzero powers of~$z$ and~$p/z$ the space of cuspidal 
$p$-adic overconvergent modular forms.
\end{proposition}
The proposition can be proved by considering what the basis of the Banach space is, and generalizing the results from the case where the genus of~$X_0(p)$ is~0.

To compute the matrix of the~$U_p$ operator with respect to the basis~\eqref{eq:orthonormal-basis}, we need to be able 
to write~$U_p(z^i)$ and~$U_p((p/z)^i)$ in terms of~$z^i$ and~$(p/z)^i$. 
Because we are dealing with modular curves with 
\emph{two} cusps in this paper, we will not be able to write these as power series in one variable with a zero at that 
cusp; indeed, we note that both~$z$ and~$p/z$ have a nonzero constant term in their $q$-expansions. This means that we 
are going to have to 
approximate the action of the~$U_p$ operator modulo~$p^N$ for higher and higher powers of~$p$ to find the 
characteristic power series of~$U_p$.

Concretely, we can find approximations to~$U_p(z^i)$, for~$i \in \Z$, by 
finding a linear combination of~$z^j$ and~$(p/z)^k$ (for~$j,k \in \N$) 
which is congruent to~$U_p(z^i)$ modulo~$p^N$ for some positive integer~$N$ using a 
computer algebra system such as~\Magma{}~\cite{magma}. This requires linear algebra over~$\Z/p^N\Z$, but is not otherwise technically difficult.

Similarly, we note that for primes~$p$ where the modular curve~$X_0(p)$ has higher genus, we would need to use more parameters to obtain a basis for the overconvergent modular forms, one for each cusp.
\section{Some computational results}

We will illustrate this with an example. Let~$p=11$ and~$k=0$. We can compute the characteristic power series of~$U_p$ 
acting on 
overconvergent modular forms of weight~0 using the procedure we have outlined above; for instance, 
if we compute the action of~$U_{11}$ modulo~$11^{13}$ acting on the basis~$\{1,z,11/z,\ldots,z^6,(11/z)^6\}$, 
we find that the characteristic power series modulo~$11^{13}$ is given by
\begin{eqnarray*}
f(x)&=&1 + 30120372860126x + 17601733022753x^2\\ &+& 32271675221764x^3 +
    17634685093520x^4 + 5939670233629x^5
\end{eqnarray*}
and it can be easily checked that the slopes of this polynomial which correspond to cusp forms are~$\{1,2,3,4\}$.
%

There are conjectures of Buzzard (see~\cite{buzzard-questions}, Section~2) which predict what the slopes of the Hecke 
operators acting on classical modular forms of weight greater than~1 will be; we can use the work of Wan~\cite{wan} on 
families of modular forms, for instance, which proves a weak version of the \Gouvea{}-Mazur conjecture (see~\cite{gouvea-mazur-conjecture}) which allows us 
to relate these slopes of weight~$k$ classical forms to the slopes of overconvergent forms of weight~0 that we compute. 
One can therefore check explicitly what the first few slopes should be in weight~0, and we find that the first few 
slopes are those we have computed above.

\subsection{Extensions to more general weights}

We have now proved that there is an orthonormal basis for the space of overconvergent modular forms of weight~0, given by~\eqref{eq:orthonormal-basis}. We can use the presentation of the nontrivial weight forms given in~\eqref{eq:weight-k-overconvergent-forms} to write down the following basis of the weight-character~$\kk$ forms:
\begin{equation}
\label{eq:orthonormal-basis-nonzero-weight}
\{E_\kk,E_\kk\cdot z,E_\kk\cdot z^2,\ldots\} \cup \{E_\kk\cdot p/z,E_\kk\cdot (p/z)^2,\ldots\}.
\end{equation}
We then repeat the same procedure to approximate the action of the~$U_p$ operator as in the previous section, 
using this basis to approximate~$U_p(E_\kk \cdot z^i)$. We see that this means that we have a more unified treatment of forms of general weight than was possible, for instance, in~\cite{kilford-2slopes}, where we had to use ``Coleman's trick'' to move between weight~0 and weight-character~$\kk$; see for instance Definition~18 of~\cite{kilford-2slopes} for an example of this.

As in the weight~0 case, we note that we can use our calculations to verify the conjectures of Buzzard and Clay on the slopes of modular forms in specific cases. We can also use this process to find the slopes of \emph{classical} modular forms, using the following well-known theorem of Coleman:
\begin{theorem}[Theorem~1.1, \cite{coleman-overconvergent}]
Let~$f$ be a $p$-adic overconvergent modular form of weight~$k$. If the slope of~$f$ is strictly less than~$k-1$, then~$f$ is a classical modular form.
\end{theorem}
Let us illustrate this with another explicit example. Let~$p=17$ and~$k=4$. We can compute the characteristic power 
series of~$U_p$ acting on overconvergent modular forms of weight~4 using the procedure above, and if we compute it 
modulo~$17^{13}$ acting on the basis~$\{E_4,E_4 \cdot z,E_4 \cdot 
17/z,\ldots,E_4 \cdot z^6,E_4 \cdot (17/z)^6\}$ then we 
find that the characteristic power series 
modulo~$17^{13}$ is given by
\begin{eqnarray*}
g(x)&=& 1 + 8750351632484700x + 8304425558400239x^2 \\ &+& 2202146495650844x^3 +
    1782357921875155x^4 + 1742180286550466x^5 \\ &+& 7082953495257174x^6 +
    6100397542758674x^7
\end{eqnarray*}
and one can check explicitly that the slopes of~$g(x)$ that correspond to cusp forms are~$\{1,1,1,1,3,4\}$; these are the same slopes as those predicted 
by 
Buzzard's conjecture. We also note that the space~$S_{4+16\cdot17}(\sltz)$ of classical modular forms is effectively computable by a computer algebra system such as \Magma{} or \Sage{}, and we find that the first few slopes of the Hecke operator~$U_{17}$ acting on these modular forms are~$\{1,1,1,1,3,4\}$; it can be seen that the largest of these classical slopes is exactly~$k-1$. 
%

We see that the methods given in this paper to compute overconvergent modular forms also enable us to draw conclusions 
about classical modular forms; this is completely independent of the modular symbols algorithms used in \Magma{}. For instance, we see from our computation above that there is a unique 17-adic overconvergent modular form of weight~4 and slope~3; this form is a classical modular form, and the uniqueness of the slope means that its Fourier expansion is defined over~$\Q_{17}$ rather than an extension field.
\subsection{Finding eigenfunctions}
A simple way to find approximations to eigenfunctions for the~$U_p$ operator by iterating the action of the matrix of~$U_p$ we have computed on a randomly-chosen basis vector, which will leave us with only the eigenforms of lowest slope.

This method is similar to that found in~\cite{gouvea-mazur-searching}, where the authors find overconvergent 5-adic modular eigenforms of weight~0 by iterating the action of the~$U_5$ operator. This in turn builds on the work of Atkin and O'Brien~\cite{atkin-obrien} which pioneered this technique for finding $p$-adic eigenforms for~$p=13$; we note that for both of these primes~$p$ $X_0(p)$ has genus~0. Another way of finding eigenfunctions is explored in Sections~4 and~5 of~\cite{loeffler}.

Let us illustrate the method with an example. Let~$p=19$ and~$k=0$. Using 
the methods we have developed above, it can be seen that the first few 
slopes corresponding to cusp forms are~$\{0, 1, 1, 2, 3, 3\}$; in particular, there is a unique lowest slope. We can compute~$\left(U_{19}\right)^9(z)$ modulo~$19^{9}$ to find the slope~0 eigenvector (we will need to exclude the other slope~0 eigenvector which is not a cusp form, but this can be done), and we find that the following combination of basis elements is an eigenvector of~$U_{19}$ modulo~$19^9$:
\begin{eqnarray*}
h(z)&=&1+30128692711z+175491384265\cdot19/z+274929953659z^2\\&+&293026132171\cdot (19/z)^2+ 95628143705z^3+30772040349\cdot (19/z)^3\\&+&77897422935z^4 +O(19^9)\cdot (19/z)^4,
\end{eqnarray*}
with eigenvalue~$158723425855$, which is indeed a unit in~$\Z/19^9\Z$.

\subsection{Compatibility with conjectures of Clay}

In Chapter~3 of Clay's thesis~\cite{clay-thesis}, conjectural formulae for the slopes of overconvergent modular forms are presented, in terms of the valuations of certain combinations of factorials. For example, the slopes of the~$U_5$ operator acting on weight~0 overconvergent forms are conjecturally given by
\begin{equation}
\label{slope-formula-5}
v_5\left(5^i\frac{(3i-1)!(3i)!}{i!(i-1)!}\right)\text{ for }i \in \N.
\end{equation}
In all cases where Buzzard's and Clay's conjectures have been checked, they give the same predictions for the slopes, and agree with the experimental evidence given by computation. We will now explain more about the basis of Clay's conjectures, and say what they are in this context.

Clay's conjectures say that in level~1 that the slopes of the~$U_p$ operator acting on overconvergent modular forms given by (the convex hull of)~$(p^2-1)/24$ formulae of a similar type to~\eqref{slope-formula-5}; in particular, for~$p=5$ this means that there is only one formula. For the simplest case we are dealing with, $p=11$, this means that there are~5 formulae of this type. We will now give these formulae explicitly; the calculations following are based upon Section~3.3 of~\cite{clay-thesis}. 

The following formulae can all be generalized to other weights at the expense of making the formulae more complicated; we present only the weight~0 formulae for the sake of simplicity.

We can write a positive integer~$i$ uniquely as~$5j+k$, where~$j$ and~$k$ are both non-negative integers and~$k$ is less than~5 (the integer~5 here is~$(11^2-1)/24$). The $i^{\rm th}$ slope is given by~$s_k(j)$;
\begin{eqnarray*}
s_k(j)=\left\{\begin{array}{ll}v_5\left(11^{4j}\frac{(6j)!(6j-1)!}{j!(j-1)!}\right), &\text{ if }k=0,\\
v_{11}\left(11^{4j-4}\frac{(6j-6)!(6j-6)!}{(j-1)!(j-1)!}\right), &\text{ if }k=1,\\
v_{11}\left(11^{4j-3}\frac{(6j-4)!(6j-5)!}{(j-1)!(j-1)!}\right), &\text{ if }k=2,\\
v_{11}\left(11^{4j-2}\frac{(6j-3)!(6j-4)!}{(j-1)!(j-1)!}\right), &\text{ if }k=3,\\
v_{11}\left(11^{4j-1}\frac{(6j-2)!(6j-3)!}{(j-1)!(j-1)!}\right), &\text{ if }k=4.
\end{array}\right.
\end{eqnarray*}
We note that no convex hulling procedure is necessary here; the formulae give the slopes without having to adjust them.

Similarly, we can write all positive integers~$i$ uniquely as $12i+j$, where~$i$ and~$j$ are again both non-negative integers. The conjectural formulae for the~$i^{\rm th}$ slope~$s_k(j)$ when~$p=17$ is given by
\[
s_k(j)=
\begin{cases}
v_{17}\left(17^{7j}\frac{(9j)!(9j-1)!}{j!(j-1)!}\right), & \text{ if }k=0,\\
v_{17}\left(17^{7j-7}\frac{(9j-9)!(9j-9)!}{(j-1)!(j-1)!}\right), & \text{ if }k=1,\\
v_{17}\left(17^{7j-9+k}\frac{(9j-11+k)!(9j-10+k)!}{(j-1)!(j-1)!}\right), & \text{ if }2 \le k \le 4,\\
v_{17}\left(17^{7j-10+k}\frac{(9j-12+k)!(9j-11+k)!}{(j-1)!(j-1)!}\right), & \text{ if }5 \le k \le 8,\\
v_{17}\left(17^{7j-11+k}\frac{(9j-13+k)!(9j-12+k)!}{(j-1)!(j-1)!}\right), & \text{ if }9 \le k \le 11,
\end{cases}
\]
and again no convex hulling is needed here.

Finally we present the conjectural formulae for~$p=19$.  We can write all positive integers~$i$ uniquely as $15i+j$, where~$i$ and~$j$ are again both non-negative integers. The conjectural formulae for the~$i^{\rm th}$ slope~$s_k(j)$ in this case is
\[
s_k(j)=
\begin{cases}
v_{19}\left(19^{8j}\frac{(10j)!(10i-1)!}{i!(i-1)!}\right), & \text{ if }k=0,\\
v_{19}\left(19^{8j-8}\frac{(10j-10)!(10j-10)!}{(j-1)!(j-1)!}\right), & \text{ if }k=1,\\
v_{19}\left(19^{8j-10+k}\frac{(10j-11+k)!(10j-12+k)!}{(j-1)!(j-1)!}\right), & \text{ if }2 \le k \le 4,\\
v_{19}\left(19^{8j-11+k}\frac{(10j-12+k)!(10j-13+k)!}{(j-1)!(j-1)!}\right), & \text{ if }5 \le k \le 7,\\
v_{19}\left(19^{8j-12+k}\frac{(10j-13+k)!(10j-14+k)!}{(j-1)!(j-1)!}\right), & \text{ if }8 \le k \le 10,\\
v_{19}\left(19^{8j-13+k}\frac{(10j-14+k)!(10j-15+k)!}{(j-1)!(j-1)!}\right), & \text{ if }11 \le k \le 12,\\
v_{19}\left(19^{8j-14+k}\frac{(10j-15+k)!(10j-16+k)!}{(j-1)!(j-1)!}\right), & \text{ if }k=13,14,
\end{cases}
\]
and as in the previous cases no convex hulling is needed.

It is interesting to note that the slope formulae for $p\in \{11,17,19\}$, where the modular curve~$X_0(p)$ has genus~1, and the slope formulae for~$p=7$, where~$X_0(p)$ has genus~0, are of the same type (there are two formulae that give the slopes for~$p=7$). This means that, viewed from this standpoint, the behaviour of the slopes is more similar than the different methods of computation needed in the genus~0 and genus~1 cases would suggest.

Based on the computations of Clay and the experimental evidence of our own computations, we make the following conjecture, which is a special case of Clay's conjectures:
\begin{conjecture}
Let~$p \in \{11,17,19\}$. Then the slopes of the~$U_p$ operator acting on overconvergent  $p$-adic modular forms of weight~0 are given by the functions~$s_k(j)$ presented above.
\end{conjecture}

\section{Acknowledgements}
I would like to thank Kevin Buzzard for suggesting this problem to me, and for helpful conversations.

I would also like to thank William Stein for giving me an account on his computer {\sc Meccah} to perform these computations.
\bibliographystyle{plain}

\begin{thebibliography}{10}

\bibitem{atkin-obrien}
A.~O.~L. Atkin and J.~N. O'Brien.
\newblock Some properties of {$p(n)$} and {$c(n)$} modulo powers of {$13$}.
\newblock {\em Trans. Amer. Math. Soc.}, 126:442--459, 1967.

\bibitem{magma}
W.~Bosma, J.~Cannon, and C.~Playoust.
\newblock The {M}agma algebra system {I}: The user language.
\newblock {\em J. Symb. Comp.}, 24(3--4):235--265, 1997.
\newblock Available from \url{http://magma.maths.usyd.edu.au}.

\bibitem{buzzard-questions}
Kevin Buzzard.
\newblock Questions about slopes of modular forms.
\newblock {\em Ast\'erisque}, 298:1--15, 2005.

\bibitem{buzzard-calegari}
Kevin Buzzard and Frank Calegari.
\newblock Slopes of overconvergent 2-adic modular forms.
\newblock {\em Compos. Math.}, 141(3):591--604, 2005.

\bibitem{buzzard-kilford}
Kevin Buzzard and L.~J.~P. Kilford.
\newblock The 2-adic eigencurve at the boundary of weight space.
\newblock {\em Compos. Math.}, 141(3):605--619, 2005.

\bibitem{clay-thesis}
Lisa Clay.
\newblock {\em Some {C}onjectures {A}bout the {S}lopes of {M}odular {F}orms}.
\newblock PhD thesis, Northwestern University, 2005.

\bibitem{coleman-overconvergent}
R.~Coleman.
\newblock Classical and overconvergent modular forms of higher level.
\newblock {\em Journal de Th\'eorie des Nombres de Bordeaux}, 9(2):395--403,
  1997.

\bibitem{coleman}
R.~Coleman.
\newblock $p$-adic {B}anach spaces and families of modular forms.
\newblock {\em Inv. Math}, 127:417--479, 1997.

\bibitem{coleman-stevens-teitelbaum}
R.~Coleman, G.~Stevens, and J.~Teitelbaum.
\newblock Numerical experiments on families of $p$-adic modular forms.
\newblock {\em AMS/IP Studies in Advanced Mathematics}, 7:143--158, 1998.

\bibitem{coleman-stein}
Robert~F. Coleman and William~A. Stein.
\newblock Approximation of eigenforms of infinite slope by eigenforms of finite
  slope.
\newblock In {\em Geometric aspects of Dwork theory. Vol. I, II}, pages
  437--449. Walter de Gruyter GmbH \& Co. KG, Berlin, 2004.

\bibitem{emerton-thesis}
M.~Emerton.
\newblock {\em 2-adic Modular Forms of minimal slope}.
\newblock PhD thesis, Harvard University, 1998.

\bibitem{gouvea-mazur-conjecture}
Fernando~Q. Gouv{\^e}a and Barry Mazur.
\newblock Families of modular eigenforms.
\newblock {\em Math. Comp.} 58 (1992), no. 198, 793--805. 

\bibitem{gouvea-mazur-searching}
Fernando~Q. Gouv{\^e}a and Barry Mazur.
\newblock Searching for {$p$}-adic eigenfunctions.
\newblock {\em Math. Res. Lett.}, 2(5):515--536, 1995.

\bibitem{katz}
Nicholas~M. Katz.
\newblock {$p$}-adic properties of modular schemes and modular forms.
\newblock In {\em Modular functions of one variable, III (Proc. Internat.
  Summer School, Univ. Antwerp, Antwerp, 1972)}, pages 69--190. Lecture Notes
  in Mathematics, Vol. 350. Springer, Berlin, 1973.

\bibitem{kilford-2slopes}
L.~J.~P. Kilford.
\newblock Slopes of 2-adic overconvergent modular forms with small level.
\newblock {\em Math. Res. Lett.}, 11(5-6):723--739, 2004.

\bibitem{kilford-5slopes}
L.~J.~P. Kilford.
\newblock On the slopes of the {$U_5$} operator acting on overconvergent
  modular forms.
\newblock {\em Journal de Th\'eorie des Nombres de Bordeaux}, 20, no. 1, 165--182, 2008.

\bibitem{kolberg}
O.~Kolberg.
\newblock {Congruences for the coefficients of the modular invariant $j (\tau)$
  modulo powers of 2.}
\newblock {\em Arbok Univ. Bergen, Mat.-naturv. Ser.}, 16:1--9, 1961.

\bibitem{loeffler}
David Loeffler.
\newblock Spectral expansions of overconvergent modular functions.
\newblock {\em Int. Math. Res. Not. IMRN}, 16:Art. ID rnm050, 17, 2007.

\bibitem{serre}
J.-P. Serre.
\newblock Endomorphismes completements continues des espaces de {B}anach
  $p$-adique.
\newblock {\em Publ. Math. IHES}, 12:69--85, 1962.

\bibitem{smithline}
L.~Smithline.
\newblock {\em Exploring slopes of $p$-adic modular forms}.
\newblock PhD thesis, University of California at Berkeley, 2000.

\bibitem{smithline-published}
Lawren Smithline.
\newblock Compact operators with rational generation.
\newblock In {\em Number theory}, volume~36 of {\em CRM Proc. Lecture Notes},
  pages 287--294. Amer. Math. Soc., Providence, RI, 2004.

\bibitem{wan}
D.~Wan.
\newblock Dimension variation of classical and $p$-adic modular forms.
\newblock {\em Inventiones Mathematicae}, 133:449--463, 1998.

\end{thebibliography}

\end{document}